\documentclass[11pt,a4paper]{article}
\usepackage{a4wide}
\usepackage{latexsym,stmaryrd}
\usepackage{amsmath,amssymb,amsthm,bm}
\usepackage[latin1]{inputenc}
\usepackage{t1enc, lmodern}
\usepackage{graphicx}
\usepackage{algorithm, algorithmicx, algpseudocode}
\usepackage[nonamebreak, round]{natbib}
\usepackage[bookmarks=false,colorlinks=false,linkcolor=blue]{hyperref}

\graphicspath{{figs/}}

\def\d{\partial}
\def\p{\varphi}
\def\E{{\mathbb E}}

\def\N{{\mathbb N}}
\def\P{{\mathbb P}}

\def\R{{\mathbb R}}

\def\ba{\bm{\alpha}}

\def\bA{\bm{A}}
\def\bu{\bm{u}}
\def\bc{\bm{c}}
\def\bW{\bm{W}}
\def \cG{\mathcal{G}}
\def \cF{\mathcal{F}}
\def \cS{\mathcal{S}}
\def \cD{\mathcal{D}}
\def \cP{\mathcal{P}}
\def \cL{\mathcal{L}}

\def \cA{\mathcal{A}}
\def\Vect{\mathop {\mathrm{Vect}}}

\def\ind#1{{\mathbf 1}_{\{#1\}}}

\def\s{\sigma}

\def\cite{~\citet}
\def \ens{[0,T]\times \R^d}
\def\keywords#1{{\bf Keywords~:~}#1}
\numberwithin{equation}{section}

\theoremstyle{plain}
\newtheorem{thm}{Theorem}

\newtheorem{remark}[thm]{Remark}
\newtheorem{definition}[thm]{Definition}

\newtheoremstyle{hypo@style}{\topsep}{15pt}{\itshape}{0pt}{\bfseries}{}{5pt}
{\thmname{#1}\thmnumber{ #2}\thmnote{ {\mdseries (#3)}}}
\theoremstyle{hypo@style}
\newtheorem{hypo}{Hypothesis}
\newtheorem{algo}{Algorithm}

\newtheoremstyle{multhypostyle}{\topsep}{\topsep}%
{\itshape}%
{}%
{\bfseries}%
{.}%
{\newline}%
{\thmname{#1}\thmnumber{ #2}\thmnote{ #3}}%

\theoremstyle{multhypostyle}

\emergencystretch=\hsize
\begin{document}

\title{A Parallel Algorithm for solving BSDEs - Application to the pricing and
hedging of American options}

\author{Céline Labart \footnote{Email:celine.labart@univ-savoie.fr}\\
Laboratoire de Mathématiques, CNRS UMR 5127\\
Université de Savoie, Campus Scientifique\\
73376 Le Bourget du Lac, France\vspace*{1cm}\\
Jérôme Lelong \footnote{Email:jerome.lelong@imag.fr}\\
Laboratoire Jean Kuntzmann,  \\
Université de Grenoble and CNRS,\\
BP 53, 38041 Grenoble, Cedex 09, France}

\date{\today}

\maketitle

\begin{abstract}  We present a parallel algorithm for solving backward
  stochastic differential equations (BSDEs in short) which are very useful
  theoretic tools to deal with many financial problems ranging from option
  pricing option to risk management. Our algorithm based on
  \cite{gobet_labart_10} exploits the link between BSDEs and non linear partial
  differential equations (PDEs in short) and hence enables to solve high
  dimensional non linear PDEs.  In this work, we apply it to the pricing and
  hedging of American options in high dimensional local volatility models, which
  remains very computationally demanding.  We have tested our algorithm up to
  dimension 10 on a cluster of 512 CPUs and we obtained linear speedups which
  proves the scalability of our implementation. \\
  
\end{abstract}

\section{Introduction}

Pricing and hedging American options with a large number of underlying assets is
a challenging financial issue. On a single processor system, it can require
several hours of computation in high dimensions. Recent advances in parallel
computing hardware such as multi--core processors, clusters and GPUs are then of
high interest for the finance community. For a couple of years, some research
teams have been tackling the parallelization of numerical algorithms for option
pricing.  \cite{thulasiram_bondarenko_02} developed a parallel algorithm using
MPI for pricing a class of multidimensional financial derivatives using a
binomial lattice approach.  \cite{huang_thulasiram_05} presented algorithms for
pricing American style Asian options using a binomial tree method.  Concerning
the parallelization of Monte--Carlo methods for pricing multi--dimensional
Bermudan/American options, the literature is quite rare. We refer to
\cite{toke_girard_06} and to \cite{dbggs_10}. Both papers propose a
parallelization through grid computing of the \cite{ibanez_04} algorithm, which
computes the optimal exercise boundary. A GPGPU approach based on quantization
techniques has recently been developed by \cite{pages_gpgpu}.\\

Our approach is based on solving Backward Stochastic Differential Equations
(BSDEs in short). As explained in the seminal paper by \cite{karoui_97}, pricing
and hedging European options in local volatility models boil down to solving
standard BSDEs. From \cite{karoui_2_97}, we know that the price of an American
option is also linked to a particular class of BSDEs called reflected BSDEs.
Several sequential algorithms to solve BSDEs can be found in the literature.
\cite{ma_94} adopted a PDE approach, whereas \cite{bouc:touz:04} and
\cite{lemor_05} used a Monte--Carlo approach based on the dynamic programming
equation. The Monte--Carlo approach was also investigated by \cite{bally_03} and
\cite{delarue_05} who applied quantization techniques to solve the dynamic
programming equation. Our approach is based on the algorithm developed by
\cite{gobet_labart_10} which combines Picard's iterations and an adaptive
control variate.  It enables to solve standard BSDEs, ie, to get the price and
delta of European options in a local volatility model. Compared to the
algorithms based on the dynamic programming equation, ours provides regular
solutions in time and space (which is coherent with the regularity of the option
price and delta).  To apply it to the pricing and hedging of American options,
we use a technique introduced by \cite{karoui_2_97}, which consists in
approximating a reflected BSDE by a sequence of standard BSDEs with
penalisation.

The paper is organized as follows. In section \ref{sec:intro}, we briefly recall
the link between BSDEs and PDEs which is the heart of our algorithm.  In section
\ref{sec:algo}, we describe the algorithm and in Section \ref{sec:parallel} we
explain how the parallelization has been carried out.  Section~\ref{sec:pricing}
describes how American options can be priced using BSDEs.  Finally, in
Section~\ref{sec:num_results}, we conclude the paper by some numerical tests of
our parallel algorithm for pricing and hedging European and American basket
options in dimension up to $10$.

\subsection{Definitions and Notations}\label{sec:def}
\begin{itemize}
\item Let $C^{k,l}_{b}$ be the set of continuously differentiable functions
  $\phi:(t,x)\in [0,T] \times \R^d$ with continuous and uniformly bounded
  derivatives w.r.t. $t$ (resp. w.r.t. $x$) up to order $k$ (resp. up to order
  $l$).
\item $C^k_p$ denotes the set of $C^{k-1}$ functions with piecewise continuous
  $k^{th}$ order  derivative.
\item For $\alpha \in ]0,1]$, $C^{k+\alpha}$ is the set of $C^k$ functions whose
  $k^{th}$ order derivative is Hölder continuous with order $\alpha$.
\end{itemize}

\section{BSDEs}\label{sec:intro}

\subsection{General results on standard BSDEs}\label{sec:BSDE}
Let $(\Omega, \mathcal{F},\mathbb{P})$ be a given probability space on which
is defined a $q$-dimensional standard Brownian motion $W$, whose natural
filtration, augmented with $\mathbb{P}$-null sets, is denoted
$(\mathcal{F}_t)_{0\le t\le T}$ ($T$ is a fixed terminal time). We denote
$(Y,Z)$ the solution of the following backward stochastic differential
equation (BSDE) with fixed terminal time $T$
\begin{equation} \label{eqs1} -dY_t=f(t,X_t,Y_t,Z_t)dt-Z_tdW_t, \;\;
  Y_T=\Phi(X_T),
\end{equation}
where $f:[0,T] \times \R^d \times \R \times \R^q \rightarrow \R$,
$\Phi:\R^d  \rightarrow \R$ and $X$ is the $\mathbb{R}^d$-valued
process solution of
\begin{equation} \label{eqs2}
  X_t= x + \int_0^t b(s,X_s) ds + \int_0^t \sigma(s,X_s)dW_s,
\end{equation}
with $b:[0,T] \times \R^d \rightarrow \R^d$ and $\sigma : [0,T] \times \R^d
\rightarrow \R^{d\times q}.$

From now on, we assume the following Hypothesis, which ensures the existence and
uniqueness of the solution to Equations~(\ref{eqs1})-(\ref{eqs2}).
\begin{hypo}\label{hypo1}\hfill
  \begin{itemize}
  \item The driver $f$ is a bounded Lipschitz continuous function, ie,  for
    all $\linebreak[4](t_1,x_1,y_1,z_1)$, $(t_2,x_2,y_2,z_2)\in [0,T] \times
    \R^d \times \R \times \R^{},\; \exists \, L_f>0, $
    \begin{align*}
      |f(t_1,x_1,y_1,z_1)-f(t_2,x_2,y_2,z_2)| \le
      L_f(|t_1-t_2|+|x_1-x_2|+|y_1-y_2|+|z_1-z_2|).
    \end{align*}
  \item $\sigma$ is uniformly elliptic on $[0,T] \times \R^d$, ie, there exist two
    positive constants $\s_0, \s_1$ s.t. for any $\xi \in \R^d$ and any $(t,x) \in
    \ens$
    \begin{align*}
      \s_0 |\xi|^2 \le \sum_{i,j=1}^d [\s \s^*]_{i,j}(t,x) \xi_i \xi_j \le \s_1
      |\xi|^2.
    \end{align*}
  \item $\Phi$ is bounded in $C^{2+\alpha}$, $\alpha \in ]0,1]$.
  \item $b$ and $\sigma$ are in $C^{1,3}_b$ and $\d_t \sigma$ is in $C^{0,1}_b$.
  \end{itemize}
\end{hypo}

\subsection{Link with semilinear PDEs}
Let us also recall the link between BSDEs and semilinear PDEs. Although the relation is
the keystone of our algorithm, as explained in Section \ref{sec:algo}, we do not
develop it and refer to \cite{pardoux_92} or \cite{karoui_97} for more details.

According to \cite[Theorem 3.1]{pardoux_92}, we can link 
the solution $(Y,Z)$ of the BSDE (\ref{eqs1}) to the solution $u$ of the
following PDE:
\begin{equation}\label{edp}
  \begin{cases}
    \d_t u(t,x) + \cL u(t,x) + f(t,x, u(t,x), (\d_x u \sigma)(t,x))=0,\\
    u(T,x)=\Phi(x),
  \end{cases}     
\end{equation}
where $\cL$ is defined by
\begin{equation*}
  \cL_{(t,x)}u(t,x)=\frac{1}{2}\sum_{i,j}[\sigma \sigma^*]_{ij}(t,x)
  \partial_{x_ix_j}^2u(t,x) +\sum_i b_i(t,x)\partial_{x_i}u(t,x).
\end{equation*}

\begin{thm}[\cite{delarue_05}, Theorem 2.1]\label{thm:dm}
  Under Hypothesis \ref{hypo1}, the solution $u$ of
  PDE (\ref{edp}) belongs to $C^{1,2}_b$. Moreover, the solution $(Y_t, Z_t)_{0 \le t \le
  T}$ of~(\ref{eqs1}) satisfies
  \begin{align}\label{eqs5}
    \forall t \in [0,T],\;\; (Y_t, Z_t)=(u(t,X_t),\d_x u(t,X_t) \sigma(t,X_t)).
  \end{align}
\end{thm}

\section{Presentation of the Algorithm}\label{sec:algo}

\subsection{Description} We present the algorithm introduced by
\cite{gobet_labart_10} to solve standard BSDEs. It is based on Picard's iterations
combined with an adaptive Monte--Carlo method. We recall that we aim at numerically
solving BSDE~(\ref{eqs1}), which is equivalent to solving the semilinear
PDE~(\ref{edp}). The current algorithm provides an approximation of the solution of
this PDE. Let $u^k$ denote the approximation of the solution $u$ of (\ref{edp}) at
step $k$. If we are able to compute an explicit solution of (\ref{eqs2}), the
approximation of $(Y,Z)$ at step $k$ follows from (\ref{eqs5}): $(Y^k_t,
Z^k_t)=(u^k(t,X_t),\partial_x u^k(t,X_t) \sigma(t,X_t))$, for all $t \in [0,T]$.
Otherwise, we introduce $X^N$ the approximation of $X$ obtained with a $N$--time step
Euler scheme:
\begin{align}\label{eqs3}
  \forall s \in [0,T],\;\;  dX^N_s=b(\p^N(s),X^N_{\p^N(s)})ds+\s(\p^N(s),X^N_{\p^N(s)})dW_s,
\end{align}
where $\p^N(s)=\sup\{t_j\,:\,t_j \le s\}$ is the largest discretization time not
greater than $s$ and $\{0=t_0 < t_1 < \cdots < t_N=T\}$ is a regular subdivision of
the interval $[0,T]$. Then, we write
\begin{align*}
  (Y^k_t, Z^k_t)=(u^k(t,X^N_t),\partial_x u^k(t,X^N_t) \sigma(t,X^N_t)),\mbox{ for
  all }t \in [0,T].
\end{align*}
It remains to explain how to build the approximation $(u^k)_k$ of $u$. The
basic idea is the following:
$$u^{k+1}=u^k+\mbox{ Monte--Carlo evaluations of the error}(u-u^k).$$

Combining Itô's formula applied to $u(s,X_s)$ and $u^k(s,X^N_s)$ between
$t$ and $T$ and the semilinear PDE (\ref{edp}) satisfied by $u$, we get that
the correction term $c^k$ is given by
\begin{align*}
  c^k(t,x)=(u-u^k)(t,x)= \mathbb{E} \left[\Psi\left(t,x,f_u,\Phi,W\right)-
  \Psi^N\left(t,x,-(\partial_t+\mathcal{L}^N)u^k,u^k(T,.),W\right)
  |\mathcal{G}^{k}\right]
\end{align*}
where
\begin{itemize}
\item $ \cL^Nu(s,X^N_s)=\frac{1}{2}\sum_{i,j} [\sigma \sigma^*]_{ij}
  (\p(s),X^N_{\p(s)})\partial_{x_ix_j}^2u(s,X^N_s)+
  \sum_ib_i(\p(s),X^N_{\p(s)})\partial_{x_i}u(s,X^N_s)$.
\item $f_v:[0,T]
  \times \R^d \rightarrow \R$ denotes the following function
  \begin{align*}
    f_v(t,x)=f(t,x,v(t,x),(\d_x v \sigma)(t,x)),
  \end{align*} where $f$ is the driver of BSDE (\ref{eqs1}), $\sigma$
  is the diffusion coefficient of the SDE satisfied by $X$ and $v:[0,T]
  \times \R^d \rightarrow \R$ is $C^1$ w.r.t. to its second argument.
\item  $\Psi$ and $\Psi^N$ denote
  \begin{align*}
    &\Psi(s,y,g_1,g_2,W)=\int_s^T g_1(r,X_r^{s,y}(W))dr +g_2(X_T^{s,y}(W)),\\
    &\Psi^N(s,y,g_1,g_2,W)=\int_s^T g_1(r,X_r^{N,s,y}(W))dr +g_2(X_T^{N,s,y}(W)),
  \end{align*}
  where $X^{s,y}$ (resp. $X^{N,s,y}$) denotes the diffusion process solving
  (\ref{eqs2}) and starting from $y$ at time $s$ (resp. its approximation
  using an Euler scheme with $N$ time steps), and $W$ denotes the standard
  Brownian motion appearing in (\ref{eqs2}) and used to simulate $X^N$, as
  given in (\ref{eqs3}).
\item $\mathcal{G}^{k}$ is the $\sigma$-algebra generated by the set of all
  random variables used to build $u^{k}$. In the above formula of $c^k$, we
  compute the expectation w.r.t. the law of $X$ and $X^N$ and not w.r.t. the
  law of $u^{k}$, which is $\mathcal{G}^{k}$ measurable. (See Definition
  \ref{def:gk} for a rigorous definition of $\mathcal{G}^{k}$). 
\end{itemize}

Note that $\Psi$ and $\Psi^N$ can actually be written as expectations by
introducing a random variable $U$ uniformly distributed on $[0,1]$.
\begin{align*}
  &\Psi(s,y,g_1,g_2,W)=\E_U\left[(T-s)g_1(s+(T-s)U,X_{s+(T-s)U}^{s,y}(W))
  +g_2(X_T^{s,y}(W))\right],\\
  &\Psi^N(s,y,g_1,g_2,W)=\E_U\left[(T-s) g_1(s+(T-s)U,X_{s+(T-s)U}^{N,s,y}(W))
  +g_2(X_T^{N,s,y}(W))\right].
\end{align*}
In the following, let $\psi^N(s,y, g_1, g_2, W, U)$ denote
\begin{equation}
  \label{psiN}
  \psi^N(s,y, g_1, g_2, W, U) = (T-s) g_1(s+(T-s)U,X_{s+(T-s)U}^{N,s,y}(W))
  +g_2(X_T^{N,s,y}(W))
\end{equation}
such that $\Psi^N(s,y,g_1,g_2,W)=\E_U[\psi^N(s,y, g_1, g_2, W, U)]$.

From a practical point of view, the PDE~\eqref{edp} is solved on $[0,T] \times \cD$
where $\cD \subset \R^d$ such that $\sup_{0 \le t \le T} |X_t| \in \cD$ with a
probability very close to $1$.
\begin{algo}\label{algo2}
  We begin with $u^0 \equiv 0$. 
  Assume that an approximated solution $u^k$ of class $C^{1,2}$ is built at step $k-1$.
  Here are the different steps to compute $u^{k+1}$.
  \begin{itemize}
  \item Pick at random $n$ points $(t^k_i,x^k_i)_{1\le i \le n}$ uniformly
    distributed over  $[0,T]\times \cD$.
  \item Evaluate  the Monte--Carlo correction $c^k$ at step $k$  at the points
    $(t^k_i,x^k_i)_{1 \le i \le n}$ using $M$ independent simulations
    \begin{align*}
      c^k(t^k_i,x^k_i)= \frac{1}{M}\sum_{m=1}^M
      \left[\psi^N\left(t^k_i,x^k_i,f_{u_{k}}+(\partial_t+\mathcal{L}^N)u_{k},
      \Phi-u^k,W^{m,k,i}, U^{m,k,i} \right)\right].
    \end{align*}

  \item Compute the vector $(u^k(t^k_i,x^k_i))_{1\le i \le n}$. Now, we know the
    vector $(u^k + c^k)(t^k_i,x^k_i))_{1\le i \le n}$. From these values, we
    extrapolate the function $u^{k+1} = u^k + c^k$ on $[0,T]\times \cD$.
    \begin{equation}
      \label{eq:unext}
      u^{k+1} (t,x) = \cP^k (u^k + c^k)(t,x), \quad \mbox{for }(t,x) \in
      [0,T]\times \cD,
    \end{equation}
    where $P^k$ is a deterministic operator, which only uses the values of the
    function at the points $(t^k_i,x^k_i)_{1\le i \le n}$ to approximate the function
    on the whole domain $[0,T] \times \cD$. The choice of the operator $P^k$ is
    discussed in Section~\ref{sec:operator}.
  \end{itemize}
\end{algo}

Since $c^k$ is computed using Monte-Carlo simulations instead of a true expectation,
the values $(c^k(t^k_i,x^k_i))_{1\le i \le n}$ are random variables. Therefore,
$u^{k+1}$ is a random function depending on the random variables needed to compute
$u_{k}$ and $W^{m,k,i}, U^{m,k,i}$, $1 \le m \le M, 1\le i \le n$. In view of this
comment, the $\sigma$-algebra $\cG^k$ has to be redefined to take into account this
new source of randomness.
\begin{definition}[Definition of the $\sigma$-algebra $\cG^k$]\label{def:gk} Let
  $\cG^{k+1}$ define the $\sigma$-algebra generated by the set of all random
  variables used to build $u^{k+1}$. Using (\ref{eq:unext}) yields
  \begin{align*}
    \cG^{k+1}=\cG^k \vee \s(\cA^k,\cS^k),
  \end{align*}  
  where $\cA^k$ is the set of random points used at step $k$ to build the estimator
  $\cP^k$ (see below), $\cS^k=\{W^{m,k,i},U^{m,k,i}, 1 \le m \le M, 1\le i \le n\}$,
  is the set of independent Brownian motions used to simulate the paths
  $X^{m,k,N}(x^k_i)$, and $\cG^k$ is the $\sigma$-algebra generated by the set of all
  random variables used to build $u^k$.
\end{definition}  

\subsection{Choice of the operator}
\label{sec:operator}

The most delicate part of the algorithm is how to extrapolate a function $h$ and its
derivatives when only knowing its values at $n$ points $(t_i,
x_i)_{i=1,\dots,n} \in [0,T] \times \cD$.

\subsubsection{A kernel operator} 

In the first version of Algorithm~\ref{algo2} presented in~\cite{gobet_labart_10}, a
function $h$ was extrapolated from the values computed on the grid by using a kernel
operator of the form
\begin{equation*}
  h (t,x) = \sum_{i=1}^n u(t_i, x_i) K_t(t-t_i) K_x(x-x_i),
\end{equation*} 
where $K_t$ is a one dimensional kernel whereas $K_x$ is a product of $d$ one
dimensional kernels. Hence, evaluating the function $h$ at a given point $(t,x)$
requires $O(n \times d)$ computations. 

The convergence result established by \cite[Theorem 5.1]{gobet_labart_10} is based on
the properties of the operator presented in \cite[Section 4]{gobet_labart_10}. Using
the linearity and the boundedness of the operator, they managed to prove that the
errors $\| v-\cP^k v\|$ and $\|\partial_x v -\partial_x (\cP^k v)\|$ are bounded,
which is a key step in proving the convergence of the algorithm.  At the end of
their paper, they  present an operator based on kernel estimators satisfying the
assumptions required to prove the convergence of the algorithm. 

\subsubsection{An extrapolating operator}
\label{extrapol}

The numerical properties of kernel operators are very sensitive to the choice of
their window parameters which is quite hard to tune for each new problem. Hence, we
have tried to use an other solution. Basically, we have borrowed the solution
proposed by~\cite{longstaff_01} which consists in extrapolating a function by solving
a least square problem defined by the projection of the original function on a
countable set of functions. Assume we know the values $(y_i)_{i=1,\dots,n}$ of a
function $h$ at the points $(t_i, x_i)_{i=1,\dots,n}$, the function $h$ can be
extrapolated by computing
\begin{equation}
  \label{lspb}
  \ba =   \arg\min_{\alpha\in\R^p} \sum_{i=1}^n \left| y_i
  - \sum_{l=1}^p \alpha_l B_l (t_i, x_i) \right|^2,
\end{equation}
where $(B_l)_{l=1,\dots,p}$ are some real valued functions defined on $[0,T] \times
\cD$. Once $\ba$ is computed, we set $\hat{h}(t,x) = \sum_{l=1}^p \ba_l B_l (t,
x)$. For the implementation, we have chosen the $(B_l)_{l=1,\dots,p}$ as a free
family of multivariate polynomials. For such a choice, $\hat{h}$ is known to converge
uniformly to $h$ when $p$ goes to infinity if $\cD$ is a compact set and $h$ is
continuous on $[0,T] \times \cD$. Our algorithm also requires to compute
the first and second derivatives of $h$ which are approximated by the first and
second derivatives of $\hat{h}$. Although the idea of approximating the derivatives of
a function by the derivatives of its approximation is not theoretically well
justified, it is proved to be very efficient in practice. We refer to
\cite{AmerGreeks} for an application of this principle to the computations of the
Greeks for American options. 

\paragraph{Practical determination of the vector $\ba$} In this part, we use the
notation $d'=d+1$. It is quite easy to see from Equation~\eqref{lspb} that $\ba$ is
the solution of a linear system. The value $\ba$ is a critical point of the criteria
to be minimized in Equation~\eqref{lspb} and the vector $\ba$ solves
\begin{align}
  \sum_{l=1}^p \ba_l \sum_{i=1}^n B_l(t_i, x_i) B_j(t_i, x_i) & = \sum_{i=1}^n y_i
  B_j(t_i, x_i) \quad \mbox{for $j=1,\dots,p$}\nonumber \\
  \label{alphasys}
  \bA \ba & = \sum_{i=1}^n y_i B(t_i, x_i) 
\end{align}
where the $p \times p$ matrix $\bA = (\sum_{i=1}^n B_l(t_i, x_i) B_j(t_i,
x_i))_{l,j=1,\dots,p}$ and the vector $B = (B_1, \dots, B_p)^*$. The matrix $\bA$ is
symmetric and positive definite but often ill-conditioned, so we cannot rely on the
Cholesky factorization to solve the linear system but instead we have to use some
more elaborate techniques such as a $QR$ factorization with pivoting or a singular
value decomposition approach which can better handle an almost rank deficient matrix.
In our implementation of Algorithm~\ref{algodetail}, we rely on the routine {\em
dgelsy} from Lapack~\cite{lapack}, which solves a linear system in the least square
sense by using some QR decomposition with pivoting combined with some
orthogonalization techniques.  Fortunately, the ill-conditioning of the matrix $\bA$
is not fate; we can improve the situation by centering and normalizing the
polynomials $(B_l)_l$ such that the domain $[0,T] \times \cD$ is actually mapped to
$[-1, 1]^{d'}$. This reduction improves the numerical behaviour of the chaos
decomposition by a great deal. \\

The construction of the matrix $\bA$ has a complexity of $O(p^2 n d')$.  The
computation of $\ba$ (Equation~\ref{eq:alpha}) requires to solve a linear system of
size $p \times p$ which requires $O(p^3)$ operations. The overall complexity for
computing $\ba$ is then $O(p^3 + n p^2 d'$).

\paragraph{Choice of the $(B_l)_l$.}

The function $u^k$ we want to extrapolate at each step of the algorithm is proved to
be quite regular (at least $C^{1,2}$), so using multivariate polynomials for the
$B_l$ should provide a satisfactory approximation.  Actually, we used polynomials
with $d'$ variates, which are built using tensor products of univariate polynomials
and if one wants the vector space $\Vect\{B_l, \; l=1,\dots,p\}$ to be the space of
$d'-$variate polynomials with global degree less or equal than $\eta$, then $p$ has to
be equal to the binomial coefficient $\binom{d'+\eta}{\eta}$. For instance, for $\eta=3$ and
$d'=6$ we find $p=84$. This little example shows that $p$ cannot be fixed by
specifying the maximum global degree of the polynomials $B_l$ without leading to an
explosion of the computational cost, we therefore had to find an other approach.  To
cope with the curse of dimensionality, we studied different strategies for truncating
polynomial chaos expansions. We refer the reader to Chapter 2 of \cite{blatman} for a
detailed review on the topic. From a computational point of view, we could not afford
the use of adaptive sparse polynomial families because the construction of the family
is inevitably sequential and it would have been detrimental for the speed-up of our
parallel algorithm. Therefore, we decided to use sparse polynomial chaos
approximation based on an hyperbolic set of indices as introduced by
\cite{BlatmanSudret}.

A canonical polynomial with $d'$ variates can be defined by a multi-index $\nu \in
\N^{d'}$ \textemdash\ $\nu_i$ being the degree of the polynomial with respect the
variate $i$. Truncating a polynomial chaos expansion by keeping only the
polynomials with total degree not greater than $\eta$ corresponds to the set of
multi-indices: $\{ \nu \in \N^{d'} : \sum_{i=1}^{d'} \nu_i \le \eta \}$. The idea of
hyperbolic sets of indices is to consider the pseudo $q-$norm of the multi-index
$\nu$ with $q \le 1$
\begin{equation}
  \label{eq:hyper_set}
  \left\{\nu \in \N^{d'} : \left( \sum_{i=1}^{d'} \nu_i^q  \right)^{1/q}\le \eta\right\}.
\end{equation}
Note that choosing $q=1$ gives the full family of polynomials with total degree not
greater than $\eta$. The effect of introducing this pseudo-norm is to favor low-order
interactions.

\section{Parallel approach}\label{sec:parallel}

In this part, we present a parallel version of Algorithm~\ref{algo2}, which is
far from being embarrassingly parallel as a crude Monte--Carlo algorithm. We
explain the difficulties encountered when parallelizing the algorithm and how we
solved them.

\subsection{Detailed presentation of the algorithm}

Here are the notations we use in the algorithm.
\begin{itemize}
\item $\bu^k = (u^k(t^{k}_i, x^{k}_i))_{1 \le i \le n} \quad \in \R^n$
\item $\bc^k = (c^k(t^k_i, x^k_i))_{1 \le i \le n} \quad \in \R^n$
\item $n$: number of points of the grid
\item $K_{it}$: number of iterations of the algorithm
\item $M$: number of Monte--Carlo samples
\item $N$: number of time steps used for the discretization of $X$
\item $p$: number of functions $B_l$ used in the extrapolating operator. This is not
  a parameter of the algorithm on its own as it is determined by fixing $\eta$ and $q$
  (the maximum total degree and the parameter of the hyperbolic multi-index set) but
  the parameter $p$ is of great interest when studying the complexity of the
  algorithm.
\item $(B_l)_{1 \le l \le p}$ is a family of multivariate polynomials used for
  extrapolating functions from a finite number of values.
\item $\ba^k \in \R^p$ is the vector of the weights of the chaos decomposition of
  $u^k$.
\item $d' = d+1$ is the number of variates of the polynomials $B_l$.
\end{itemize}

\begin{algorithm}[ht]
  \caption{Iterative algorithm}\label{algodetail}
  \begin{algorithmic}[1]
    \State $u^0 \equiv 0$, $\ba^0 \equiv 0$.
    \For{$k=0:K_{it}-1$} \label{Kit}
    \State Pick at random $n$ points $(t^k_i,x^k_i)_{1 \le i \le n}$.
    \For{$i=1:n$} \label{loopi}
    \For{$m=1:M$} \label{loopm}
    \State \label{Euler} Let $\bW$ be a Brownian motion with values in
    $\R^d$ discretized \newline \hspace*{2cm} on a time grid with $N$ time steps.
    \State Let $U \sim {\cal U}_{[0,1]}$.
    \State \label{MC} Compute 
    \begin{equation*}
      a_m^{i,k} = \psi^N\left(t^k_i,x^k_i,f_{u^{k}}+
      (\partial_t+\mathcal{L}^N)u^{k},\Phi-u^k,\bW, U
      \right).
    \end{equation*}
     \hspace*{2cm}$/*$ {\it We recall that $ u^{k}(t,x) = \sum_{l=1}^p \ba^{k}_l B_l(t,x)$} $*/$
  \EndFor
  \begin{align}
    \label{eq:chat}
    \bc^k_i & = \frac{1}{M}\sum_{m=1}^M a_m^{i,k} \\
    \label{eq:cu}
    \bu^k_i & = \sum_{l=1}^p \ba^{k}_l B_l(t^k_i, x^k_i)
  \end{align}
\EndFor
\State \label{nextu} Compute 
\begin{equation}
  \label{eq:alpha}
  \ba^{k+1} =   \arg\min_{\alpha \in \R^p} \sum_{i=1}^n \left| (\bu^k_i +
  \bc^k_i) - \sum_{l=1}^p \alpha_l B_l (t^k_i, x^k_i) \right|^2.
\end{equation}
  \EndFor
\end{algorithmic}
\end{algorithm}

\subsection{Complexity of the algorithm}

In this section, we study in details the different parts of
Algorithm~\ref{algodetail} to determine their complexities. Before diving into the
algorithm, we would like to briefly look at the evaluations of the
function $u^k$ and its derivatives. We recall that 
\begin{equation*}
  u^{k}(t,x) = \sum_{l=1}^p \ba^{k}_l B_l(t,x)
\end{equation*}
where the $B_l(t,x)$ are of the form $t^{\beta_{l,0}} \prod_{i=1}^d
{x_i}^{\beta_{l,i}}$ and the $\beta_{l,i}$ are some integers. Then the computational
time for the evaluation of $u^{k}(t,x)$ is proportional to $p \times d'$. The
first and second derivatives of $u^k$ write 
\begin{align*}
  \nabla_x u^{k}(t,x) & = \sum_{l=1}^p \ba^{k}_l \nabla_x B_l(t,x), \\
  \nabla_x^2 u^{k}(t,x) & = \sum_{l=1}^p \ba^{k}_l \nabla_x^2 B_l(t,x),
\end{align*}
and the evaluation of $\nabla_x B_l(t,x)$ (resp. $\nabla_x^2 B_l(t,x)$) has a
computational cost proportional to $d^2$ (resp. $d^3$). 

\begin{itemize}
\item The computation (at line 6) of the discretization of the $d-$dimensional
  Brownian motion with $N$ time steps requires $O(N d)$ computations.
\item The computation of each $a_m^{k,i}$ (line~\ref{MC}) requires the evaluation of
  the function $u^k$ and its first and second derivatives which has a cost $O(p
  d^3)$.  Then, the computation of $\bc_i^k$ for given $i$ and $k$ has a
  complexity of $O(M p d^3)$.
\item The computation of $\ba$ (Equation~\ref{eq:alpha}) requires $O(p^3 + n p^2
  d)$ operations as explained in Section~\ref{extrapol}.
\end{itemize}

The overall complexity of Algorithm~\ref{algodetail} is $O(K_{it} n M (p d^3 + d N) +
K_{it}  (p^2 n d + p^3))$. \\

To parallelize an algorithm, the first idea coming to mind is to find loops
with independent iterations which could be spread out on different processors
with a minimum of communications. In that respect, an embarrassingly parallel
example is the well-known Monte-Carlo algorithm. Unfortunately,
Algorithm~\ref{algodetail} is far from being so simple. The iterations of the
outer loop (line~\ref{Kit}) are linked from one step to the following,
consequently there is no hope parallelizing this loop. On the contrary, the
iterations over $i$ (loop line~\ref{loopi}) are independent as are the ones
over $m$ (loop line~\ref{loopm}), so we have at hand two candidates to
implement parallelizing. We could even think of a 2 stage parallelism : first
parallelizing the loop over $i$ over a small set of processors and inside this
first level parallelizing the loop over $m$. Actually, $M$ is not large enough
for the parallelization of the loop over $m$ to be efficient (see Section
\ref{sec:operator}).  It turns out to be far more efficient to parallelize
the loop over $i$ as each iteration of the loop requires a significant amount of
work.

\subsection{Description of the parallel part}

As we have just explained, we have decided to parallelize the loop over $i$
(line~\ref{loopi} in Algorithm~\ref{algodetail}). We have used a \textit{Robbin Hood}
approach. In the following, we assume to have $P+1$ processors at hand with $n>P$.
We use the following master slave scheme:
\begin{enumerate}
\item  Send to each of the $P$ slave processors the solution $\ba^k$ computed at
  the previous step of the algorithm.
\item Spread the first $P$ points $(t^k_i,x^k_i)_{1\le i \le P}$ to the $P$ slave
  processors and assign each of them the computation of the corresponding $\bc^k_i$.
\item As soon as a processor finishes its computation, it sends its result back to
  the master which in turn sends it a new point $(t^k_i, x^k_i)$ at which evaluating
  $\psi^N$ to compute $\bc^k_i$ and this process goes on until all the
  $(\bc^k_i)_{i=1,\dots,n}$ have been computed.
\end{enumerate}
At the end of this process, the master knows $\bc^{k}$, which corresponds to
the approximation of the correction at the random points $(t^k_i, x^k_i)_{1
\le i \le n}$. From these values and the vector $\bu^k$, the master computes
$\ba^{k+1}$ and the algorithm can go through a new iteration over $k$.

What is the best way to send the data to each slave process?
\begin{itemize}
\item Before starting any computations, assign to each process a block of iterations
  over $i$ and send the corresponding data all at once. This way, just one connection
  has to be initialised which is faster. But the time spent by the master to take
  care of its slave is longer which implies that at the beginning the slave process
  will remain longer unemployed. When applying such a strategy, we implicitly assume
  that all the iterations have the same computational cost.
\item Send data iteration by iteration. The latency at the beginning is smaller than
  in the block strategy and performs better when all iterations do not have the same
  computational cost.
\end{itemize}
Considering the wide range of data to be sent and the intensive use of elaborate
structures, the most natural way to pass these objects was to rely on the packing
mechanism of MPI. Moreover, the library we are using in the code (see
Section~\ref{sec:pnl}) already has a MPI binding which makes the manipulation of the
different objects all the more easy. The use of packing enabled to reduce the number
of communications between the master process and a slave process to just one
communication at the beginning of each iteration of the loop over $i$ (line~\ref{loopi} of
Algorithm~\ref{algodetail}).

\subsection{Random numbers in a parallel environment}

One of the basic problem when solving a probabilistic problem in parallel
computing is the generation of random numbers. Random number generators are
usually devised for sequential use only and special care should be taken in
parallel environments to ensure that the sequences of random numbers generated
on each processor are independent. We would like to have minimal communications
between the different random number generators, ideally after the
initialisation process, each generator should live independently of the others.

\begin{figure}[h!t]
  \centering
  \includegraphics[scale=0.5]{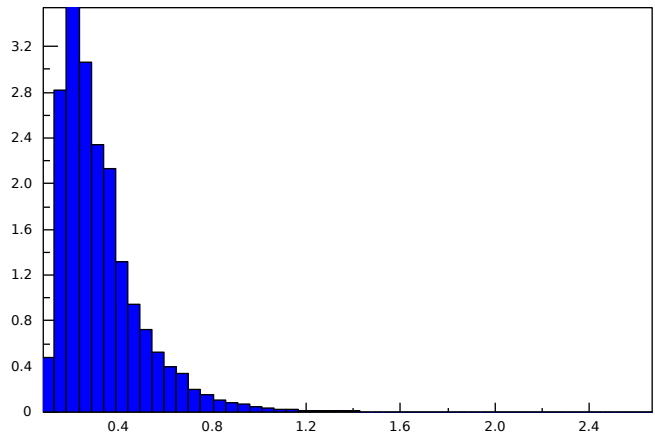}
  \caption{Distribution of the CPU time needed for the creation of one
  Mersenne Twister generator}
  \label{fig:cpuonemt}
\end{figure}

Several strategies exist for that.
\begin{enumerate}
\item Newbies in parallel computing might be tempted to take any random number
  generator and fix a different seed on each processor. Although this naive
  strategy often leads to satisfactory results on toy examples, it can induce
  an important bias and lead to detrimental results. 
\item The first reasonable approach is to split a sequence of random numbers across
  several processors. To efficiently implement this strategy, the generator must have
  some splitting facilities such that there is no need to draw all the samples prior
  to any computations. We refer to~\cite{splitting,streams} for a presentation of a
  generator with splitting facilities. To efficiently split the sequence, one should
  know in advance the number of samples needed by each processor or at least an upper
  bound of it. To encounter this problem, the splitting could be made in substreams
  by jumping ahead of $P$ steps at each call to the random procedure if $P$ is the
  number of processors involved. This way, each processor uses a sub-sequence of the
  initial random number sequence rather than a contiguous part of it. However, as
  noted by \cite{longrangecorel}, long range correlations in the original sequence
  can become short range correlations between different processors when using
  substreams. 

  Actually, the best way to implement splitting is to use a generator with a huge
  period such as the Mersenne Twister (its period is $2^{19937} -1$) and divide the
  period by a million or so if we think we will not need more than a million
  independent substreams. Doing so, we come up with substreams which still have an
  impressive length, in the case of the Mersenne Twister each substream is still
  about $2^{19917}$ long. 

\item A totally different approach is to find generators which can be easily
  parametrised and to compute sets of parameters ensuring the statistical
  independence of the related generators. Several generators offer such a
  facility such as the ones included in the SPRNG package
  (see~\cite{Mascagni97} for a detailed presentation of the generators
  implemented in this package) or the dynamically created Mersenne Twister
  (DCMT in short), see~\cite{dcmt}.
\end{enumerate}
For our experiments, we have decided to use the DCMT. This generator has a
sufficiently long period ($2^{521}$ for the version we used) and we can create at
most $2^{16} = 65536$ independent generators with this period which is definitely
enough for our needs. Moreover, the dynamic creation of the generators follows a
deterministic process (if we use the same seeds) which makes it reproducible. The
drawback of the DCMT is that its initialization process might be quite lengthy,
actually the CPU time needed to create a new generator is not bounded.  We give in
Figure~\ref{fig:cpuonemt} the distribution. 

\subsection{The library used for the implementation}
\label{sec:pnl}

Our code has been implemented in C using the PNL library (see~\cite{pnl}). This is a
scientific library available under the Lesser General Public Licence and it offers
various facilities for implementing mathematics and more recently some MPI bindings
have been added to easily manipulate the different objects available in PNL such as
vectors and  matrices. In our problem, we needed to manipulate matrices and vectors
and pass them from the master process to the slave processes and decided to use
the packing facilities offered by PNL through its MPI binding. The technical part was
not only message passing but also random number generation as we already mentioned
above and PNL offers many functions to generate random vectors or matrices using
several random number generators among which the DCMT. 

Besides message passing, the algorithm also requires many other facilities such as
multivariate polynomial chaos decomposition which is part of the library. For the
moment, three families of polynomials (Canonical, Hermite and Tchebichev polynomials)
are implemented along with very efficient mechanism to compute their first and second
derivatives. The implementation tries to make the most of code factorization to avoid
recomputing common quantities several times. The polynomial chaos decomposition
toolbox is quite flexible and offers a reduction facility such as described in
Section~\ref{extrapol} which is completely transparent from the user's side. To face
the curse of dimensionality, we used sparse polynomial families based on an
hyperbolic set of indices.

\section{Pricing and Hedging American options in local volatility models}
\label{sec:pricing}
In this Section, we present a method to price and hedge American options by using
Algorithm~\ref{algodetail}, which solves standard BSDEs.

\subsection{Framework}\label{sec:model}
Let $(\Omega,\cF ,\P)$ be a probability space and $(W_t)_{t \ge 0}$ a standard
Brownian motion with values in $\R^d$. We denote by $(\cF_t)_{t\ge 0}$ the
$\P$-completion of the natural filtration of $(W_t)_{t \ge 0}$.\\
Consider a financial market with $d$ risky assets, with prices
$S^1_t,\cdots,S^d_t$ at time $t$, and let $X_t$ be the $d$-dimensional vector
of log-returns $X^i_t=\log S^i_t$, for $i=1,\cdots,d$. We assume that
$(X_t)$ satisfies the following stochastic differential equation:
\begin{align}\label{eq4} dX^i_t=\left(r(t)-\delta_i(t)-\frac{1}{2}\sum_{j=1}^d
  \s_{ij}²(t,e^{X_t})\right)dt+\sum_{j=1}^d\sigma_{ij}(t,e^{X_t})dW^j_t, \;
  i=1,\cdots,d
\end{align}
on a finite interval $[0,T]$, where $T$ is the maturity of the option. We
denote by $X^{t,x}_s$ a continuous version of the flow of the stochastic
differential Equation (\ref{eq4}). $X^{t,x}_t=x$ almost surely.

In the following, we assume
\begin{hypo}\label{hypo2}\hfill
  \begin{enumerate}
  \item $r: [0,T] \longmapsto \R$ is a $C^{1}_b$ function. $\delta: [0,T]
    \longmapsto \R^d$ is a $C^{1}_b$ function.
  \item $\sigma : [0,T]\times \R^d \longmapsto \R^{d \times d}$ is a $C^{1,2}_b$ function.
  \item $\sigma$ satisfies the following coercivity property:
    \begin{align*}
      \exists \varepsilon>0 \; \forall (t,x) \in [0,T] \times \R^d, \forall \xi \in \R^d
      \sum_{1\le i,j\le d} [\s \s^*]_{i,j}(t,x)\xi_i\xi_j \ge \varepsilon \sum_{i=1}^d
      \xi_i^2
    \end{align*}
  \end{enumerate}
\end{hypo}
We are interested in computing the price of an American option with payoff
$\Phi(X_t)$, where $\Phi : \R^d \longmapsto \R_+$ is a continuous function (in
case of a put option, $\Phi(x)=(K-\frac{1}{d}(e^{x_1}+\cdots+e^{x_d})_+)$.

From \cite{jaillet_90}, we know that if
\begin{hypo}\label{hypo3}
  $\Phi$ is continuous and satisfies $|\Phi(x)|\le M e^{M|x|}$ for some $M>0$,
\end{hypo}
the price at time $t$ of the American option with payoff $\Phi$ is given by
\begin{align*}
  V(t,X_t)=\mbox{esssup}_{\tau \in \mathcal{T}_{t,T}}\E\left(e^{-\int_t^{\tau} r(s)
  ds} \Phi(X_{\tau})|\cF_t \right),
\end{align*} 
where $\mathcal{T}_{t,T}$ is the set of all stopping times with values in $[t,T]$ and
$V: [0,T]\times \R^d \longmapsto \R_+$ defined by $V(t,x)=\mbox{sup}_{\tau \in
\mathcal{T}_{t,T}}\E\left(e^{-\int_t^{\tau} r(s) ds} \Phi(X^{t,x}_{\tau})
\right)$.\\
We also introduce the amount $\Delta(t,X_t)$ involved in the asset at time $t$
to hedge the American option. The function $\Delta$ is given by
$\Delta(t,x)=\nabla_x V(t,x)$.

In order to link American option prices to BSDEs, we need three steps:
\begin{enumerate}
\item writing the price of an American option as a solution of a variational
  inequality (see Section \ref{sec:AO_VI})
\item linking the solution of a variational inequality to the solution of a
  reflected BSDE (see Section \ref{sec:VI_RBSDE})
\item approximating the solution of a RBSDE by a sequence of standard BSDEs
  (see Section \ref{sec:RBSDE_BSDE})
\end{enumerate}
We refer to \cite{jaillet_90} for the first step and to \cite{karoui_2_97} for
the second and third steps.

\subsection{American option and Variational Inequality}\label{sec:AO_VI}
First, we recall the variational inequality satisfied by $V$. We refer to
\cite[Theorem 3.1]{jaillet_90} for more details. Under Hypotheses
\ref{hypo2} and \ref{hypo3}, $V$ solves the following
parabolic PDE
\begin{align}\label{eqs4}
  \left\{\begin{array}{l}
    \max(\Phi(x)-u(t,x),\d_t u(t,x)+\cA u(t,x)-r(t) u(t,x))=0,\\
    u(T,x)=\Phi(x).
  \end{array}
  \right.
\end{align} where  $\cA u(t,x)=\sum_{i=1}^d
(r(t)-\delta_i(t)-\frac{1}{2}\sum_{j=1}^d
\sigma_{ij}^2(t,e^x))\d_{x_i}u(t,x)+\frac{1}{2}\sum_{1 \le i,j \le d} [\s
\s^*]_{ij}(t,e^x)\d^2_{x_ix_j} u(t,x)$ is the generator of $X$.

\subsection{Variational Inequality and Reflected BSDEs}\label{sec:VI_RBSDE}
This section is based on \cite[Theorem 8.5]{karoui_2_97}.
We assume
\begin{hypo}\label{hypo4}
  $\Phi$ is continuous and has at most a polynomial growth (ie, $\exists C,p > 0$
  s.t. $|\Phi(x)|\le C(1+|x|^p)$).
\end{hypo}
\noindent Let us consider the following reflected BSDE
\begin{align}\label{eqs6}
  \begin{cases}
    -dY_t=-r(t)Y_t dt - Z_t dW_t + dH_t,\\
    Y_T=\Phi(X_T), \; Y_t \ge \Phi(X_t), \forall t,\\
    \int_0^T (Y_t-\Phi(X_t))dH_t =0.
  \end{cases}
\end{align} Then, under Hypotheses \ref{hypo2} and \ref{hypo4},
$u(t,x)=Y_t^{t,x}$ is a viscosity solution of the obstacle problem
(\ref{eqs4}), where $Y^{t,x}_t$ is the value at time $t$ of the solution of
(\ref{eqs6}) on $[t,T]$ where the superscripts $t,x$ mean that $X$ starts from
$x$ at time $t$. We also have $(\nabla_x u(t,x))^* \sigma(t,x)=Z^{t,x}_t$ ($^*$
means transpose). Then, we get that the price $V$ and the delta $\Delta$ of
the option are given by
\begin{align*}
  V(t,X_t)=Y_t,\; \Delta(t,X_t)= (Z_t\sigma(t,X_t)^{-1})^*.
\end{align*}

\subsection{Approximation of a RBSDE by a sequence of standard
  BSDEs}\label{sec:RBSDE_BSDE}

We present a way of approximating a RBSDE by a sequence of standard BSDEs. The
idea was introduced by \cite[Section 6]{karoui_2_97} for proving the
existence of a solution to RBSDE by turning the constraint $Y_t \ge \Phi(X_t)$
into a penalisation. Let us consider the following sequence of BSDEs indexed
by $i$
\begin{align}
  \label{eq:bsde-penal}
  Y^i_t=\Phi(X_T)-\int_t^T r(s)Y^i_s ds+i\int_t^T
  (Y^i_s-\Phi(X_s))^-ds-\int_t^T Z^i_s dW_s,
\end{align} whose solutions are denoted $(Y^i,Z^i)$. We define
$H^i_t=i\int_t^T (Y^i_s-\Phi(X_s))^-ds$. Under Hypotheses \ref{hypo2} and
\ref{hypo4}, the sequence $(Y^i,Z^i,H^i)_i$ converges to the solution
$(Y,Z,H)$ of the RBSDE~(\ref{eqs6}), when $i$ goes to infinity. Moreover,
$Y^i$ converges increasingly to $Y$. The term $H^i$ is often called a
penalisation.  From a practical point, there is no use solving such a sequence
of BSDEs, because we can directly apply our algorithm to solve
Equation~\eqref{eq:bsde-penal} for a given $i$. Therefore, we actually
consider the following penalized BSDE
\begin{align}
  \label{eq:bsde-penal-2}
  Y_t=\Phi(X_T)-\int_t^T r(s)Y_s ds+\omega \int_t^T
  (Y_s-\Phi(X_s))^-ds-\int_t^T Z_s dW_s,
\end{align} where $\omega \ge 0$ is penalization weight. In practice, the
magnitude of $\omega$ must remain reasonably small as it appears as a
contraction constant when studying the speed of convergence of our
algorithm. Hence, the larger $\omega$ is, the slower our algorithm
converges. So, a trade-off has to be found between the convergence of our
algorithm to the solution $Y$ of \eqref{eq:bsde-penal-2} and the accuracy of
the approximation of the American option price by $Y$.

\subsection{European options} Let us consider a European option with payoff
$\Phi(X_T)$, where $X$ follows (\ref{eq4}). We denote by $V$ the option price and by
$\Delta$ the hedging strategy associated to the option.  From \cite{karoui_97}, we
know that the couple $(V,\Delta)$ satisfies
\begin{align*}
   -dV_t=-r(t)V_t dt - \Delta_t^* \sigma(t,X_t) dW_t, V_T=\Phi(X_T).
 \end{align*} Then, $(V,\Delta)$ is solution of a standard BSDE. This
 corresponds to the particular case $\omega=0$ of (\ref{eq:bsde-penal-2}), $Y$
 corresponding to $V$ and $Z$ to $ \Delta_t^* \sigma(t,X_t)$.

\section{Numerical results and performance}
\label{sec:num_results}

\subsection{The cluster}

All our performance tests have been carried out on a $256-$PC cluster from
SUPELEC Metz. Each node is a dual core processor : INTEL Xeon-3075 2.66 GHz
with a front side bus at 1333Mhz.  The two cores of each node share 4GB of RAM
and all the nodes are interconnected using a Gigabit Ethernet network. In none
of the experiments, did we make the most of the dual core architecture since
our code is one threaded. Hence, in our implementation a dual core processor
is actually seen as two single core processors.

The accuracy tests have been achieved using the facilities offered by the
University of Savoie computing center MUST.

\subsection{Black-Scholes' framework}

We consider a $d-$dimensional Black-Scholes model in which the dynamics
under the risk neutral measure of each asset $S^i$ is supposed to be given by
\begin{equation}\label{eq:bs}
  dS_t^i = S_t^i ((r - \delta_i)dt + \sigma^i dW^i_t) \qquad S_0 = (S_0^1, \dots, S_0^d)
\end{equation}
where $W = (W^1, \dots, W^d)$. Each component $W^i$ is a standard Brownian
motion. For the numerical experiments, the covariance structure of $W$ will be
assumed to be given by $\langle W^i, W^j \rangle_t = \rho t \ind{i \neq j} + t
\ind{i = j}$. We suppose that $\rho \in (-\frac{1}{d-1}, 1)$, which
ensures that the matrix $C=(\rho\ind{i \neq j} +\ind{i = j})_{1\leq i,j\leq
  d}$ is positive definite. Let $L$ denote the lower triangular matrix
involved in the Cholesky decomposition $C=LL^*$. To simulate $W$ on the
time-grid $0<t_1<t_2<\hdots<t_N$, we need $d\times N$ independent
standard normal variables and set
$$\begin{pmatrix}
  W_{t_1} \\ W_{t_2}\\ \vdots \\ W_{t_{N-1}} \\ W_{t_N}
\end{pmatrix}=
\begin{pmatrix}
  \sqrt{t_1}L & 0 & 0 &\hdots &0\\
  \sqrt{t_1}L &\sqrt{t_2-t_1}L & 0 &\hdots &0\\
  \vdots&\ddots&\ddots&\ddots&\vdots\\
  \vdots&\ddots&\ddots& \sqrt{t_{N-1}-t_{N-2}}L&   0 \\
  \sqrt{t_1}L & \sqrt{t_2-t_1}L  &\hdots & \sqrt{t_{N-1}-t_{N-2}}L
  &\sqrt{t_N-t_{N-1}}L
\end{pmatrix}G,$$ where $G$ is a normal random vector in $\R^{d \times N}$.
The vector $\sigma=(\sigma^1, \dots, \sigma^d)$ is the vector of volatilities,
$\delta = (\delta^1, \dots, \delta^d)$ is the vector of instantaneous dividend
rates and $r>0$ is the instantaneous interest rate. We will denote the
maturity time by $T$.  Since we know how to simulate the law of $(S_t,S_T)$
exactly for $t<T$, there is no use to discretize equation (\ref{eq:bs})
 using the Euler scheme. In this Section $N=2$.

\subsubsection{European options}

We want to study the numerical accuracy of our algorithm and to do that we first
consider the case of European basket options for which we can compute benchmark price
by using very efficient Monte-Carlo methods, see~\cite{jourdainlelong} for instance,
while it is no more the case for American options.\\
In this paragraph, the parameter $\omega$ appearing in (\ref{eq:bsde-penal-2}) is
$0$.

\paragraph{European put basket option.}
Consider the following put basket option with maturity $T$
\begin{equation}
  \label{eq:put_payoff}
  \left( K - \frac{1}{d} \sum_{i=1}^d S_T^i  \right)_+
\end{equation}

\begin{figure}[h!t]
  \centering
  \includegraphics[scale=0.6]{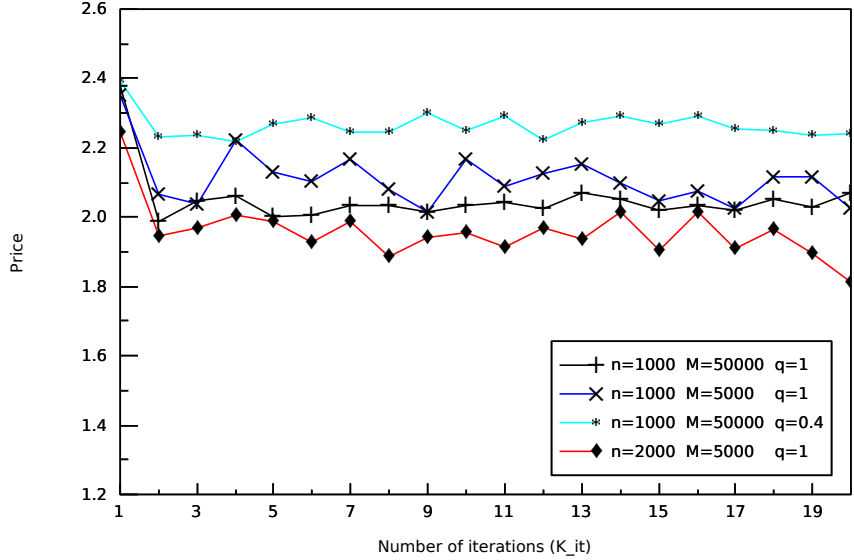}
  \caption{{\small Convergence of the algorithm for a European put basket option with $d=5$,
      $\rho=0.1$, $T=3$, $S_0 = 100$, $\sigma = 0.2$, $\delta=0$, $r=0.05$ $K=100$,
      $\eta=3$, $\omega=0$. The benchmark price computed with a high precision
      Monte--Carlo method yields $2.0353$ with a confidence interval of $(2.0323,
  2.0383)$.}}
  \label{fig:put_euro_5}
\end{figure} Figure~\ref{fig:put_euro_5} presents the influence of the
parameters $M$, $n$ and $q$. The results obtained for $n=1000$, $M=50,000$ and
$q=1$ (curve ($+$)) are very close to the true price and moreover we can see
that the algorithm stabilizes after very few iterations (less than
$10$). 

\begin{itemize}
\item Influence of $M$: curves ($+$) and ($\times$) show that taking $M=5000$ is
  not enough to get the stabilization of the algorithm.
\item Joined influence of $n$ and $M$: curves ($\times$) and ($\blacklozenge$) show
  the importance of well balancing the number of discretization points $n$ with the
  number of Monte--Carlo simulations $M$. The sharp behaviour of the curve
  ($\blacklozenge$) may look surprising at first, since we are tempted to think that
  a larger numbers of points $n$ will increase the accuracy. However, increasing the
  number of points but keeping the number of Monte--Carlo simulations constant
  creates an over fitting phenomenon because the Monte--Carlo errors arising at each
  point are too large and independent and it leads the approximation astray.
\item Influence of $q$: we can see on curves ($+$) and ($*$) that decreasing the
  hyperbolic index $q$ can lead a highly biased although smooth
  convergence. This highlights the impact of the choice of $q$ on the solution
  computed by the algorithm.
\end{itemize}

To conclude, we notice that the larger the number of Monte--Carlo simulations
is, the smoother the convergence is, but when the polynomial family considered
is too sparse it can lead to a biased convergence.

\paragraph{European call basket option.}
Consider the following put basket option with maturity $T$
\begin{equation}
  \label{eq:call_payoff}
  \left(\frac{1}{d} \sum_{i=1}^d S_T^i  - K \right)_+
\end{equation}

\begin{figure}[h!t]
  \centering
  \includegraphics[scale=0.6]{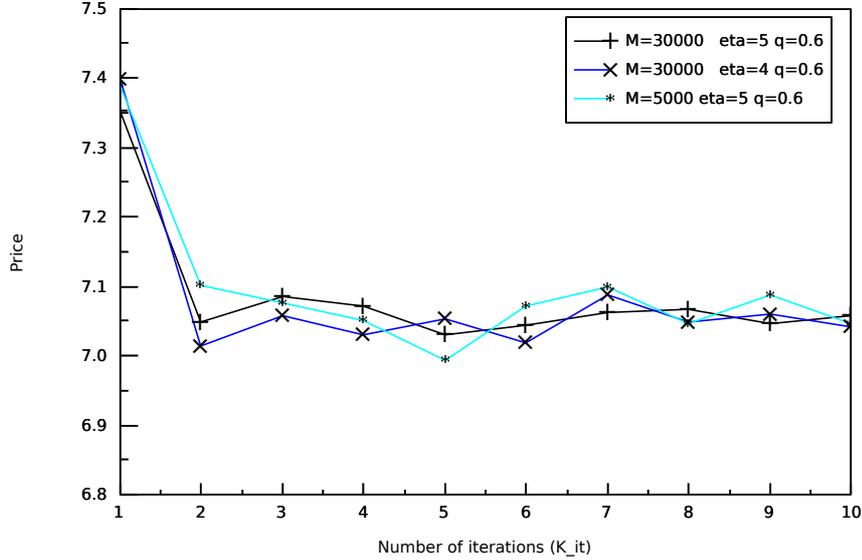}
  \caption{{\small Convergence of the price of a European call basket option with
      $d=10$, $\rho=0.2$, $T=1$, $S_0 = 100$, $\sigma = 0.2$, $\delta=0$, $r=0.05$, $K=100$,
      $n=1000$ and $\omega=0$. The benchmark price computed with a high precision
      Monte--Carlo method yields $7.0207$ with a confidence interval of $(7.0125,
  7.0288)$.}}
  \label{fig:call_euro_10}
\end{figure}

Figure~\ref{fig:call_euro_10} illustrates the impact of the sparsity of the
polynomial basis considered on the convergence of the algorithm.
 The smoothest
convergence is achieved by the curve ($+$), ie when $M=30,000$, $\eta=5$ and
$q=0.6$. The algorithm stabilizes very close to the true price and after very
few iterations.
\begin{itemize}
\item Influence of $\eta$ : for a fixed value of $q$, the sparsity increases when
  $\eta$ decreases, so the basis with $\eta=4, q=0.6$ is more sparse than the one
  with $\eta=5, q=0.6$. We compare curves ($+$) ($\eta=5$) and ($\times$) ($\eta=4$)
  for fixed values of $q$ ($=0.6$) and $M$ ($=30,000$). We can see that for $\eta=4$
  (curve ($\times$)) the algorithm stabilizes after $7$ iterations, whereas for
  $\eta=6$ (curve ($+$)) less iterations are needed to converge.
\item Influence of $M$ : for fixed values of $\eta$ $ (=5)$ and $q$ $ (=0.6)$, we
  compare curves ($+$) ($M=30000$) and ($*$) ($M=5000$). Using a large number of
  simulations is not enough to get a good convergence, as it is shown by curve ($*$).
\end{itemize}
Actually, when the polynomial basis
becomes too sparse, the approximation of the solution computed at each step of the
algorithm incorporates a significant amount a noise which has a similar effect to
reducing the number of Monte--Carlo simulations. This is precisely what we observe on
Figure~\ref{fig:call_euro_10}: the curves ($\times$) and ($*$) have a very
similar behaviour although one of them uses a much larger number of simulations.
  
\subsubsection{American options}

In this paragraph, the penalisation parameter $\omega$ appearing in
(\ref{eq:bsde-penal-2}) is $1$.

\paragraph{Pricing American put basket options.}

We have tested our algorithm on the pricing of a multidimensional American put option
with payoff given by Equation~\eqref{eq:put_payoff}

\begin{figure}[h!t]
  \centering
  \includegraphics[scale=0.5]{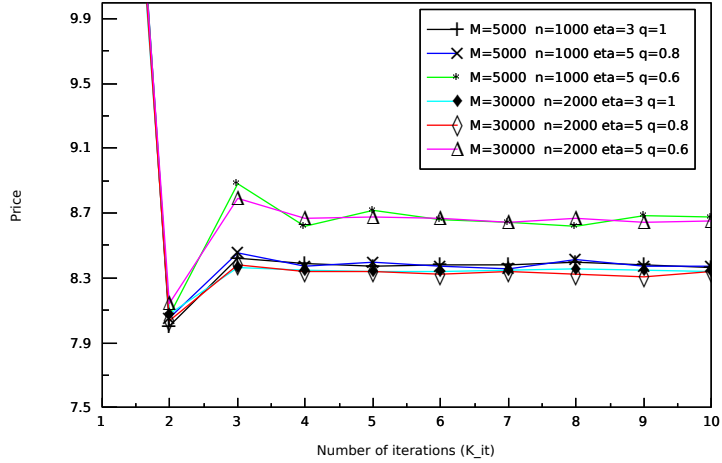}
  \caption{{\small Convergence of the price of an American put basket option with $d=5$,
  $\rho=0.2$, $T=1$, $S_0 = 100$, $\sigma = 0.25$, $\delta=0.1$, $K=100$, $r=0.05$
  and $\omega=1$.}}
    \label{fig:put_amer_5}
\end{figure}

Figure~\ref{fig:put_amer_5} presents the influence of the parameters $M$ and
$q$.
\begin{itemize}
\item Influence of $M$ : when zooming on Figure~\ref{fig:put_amer_5}, one can indeed
  see that the curves using $30000$ Monte--Carlo simulations are a little smoother
  than the others but these extra simulations do not improve the convergence as much
  as in Figures~\ref{fig:put_euro_5} and~\ref{fig:call_euro_10} (compare curves ($+$)
  and ($\blacklozenge$), Figure~\ref{fig:put_amer_5}). The main explanation of this
  fact is that put options have in general less variance than call options and in
  Figure~\ref{fig:put_euro_5} a maturity of $T=3$ was used which leads to a larger
  variance than with $T=1$.
\item Influence of $q$ : once again, we can observe that increasing the sparsity of
  the polynomial basis (ie, decreasing $q$) can lead to a biased convergence. When
  $q=0.6$, we get a biased result (see curves ($*$) and ($\vartriangle$)), even for
  $M$ large (curve ($\vartriangle$), $M=30000$).
\end{itemize}

Then, it is advisable for American put options to use almost
  full polynomial basis with fewer Monte--Carlo simulations in order to master
  the computational cost rather than doing the contrary.

\paragraph{Hedging American put basket options.}

Now, let us present the convergence of the approximation of the delta at time $0$.
Table~\ref{tab:delta_put_amer_5} presents the values of the
delta of an American put basket option when the iterations increase. We see
that the convergence is very fast (we only need $3$ iterations to get a
stabilized value). The parameters of the algorithm are the following ones:
$n=1000$, $M=5000$, $q=1$, $\eta=3$ and $\omega=1$.
\begin{table}[ht!]
  \centering
\begin{tabular}{|c|c|c|c|c|c|}
  \hline
  Iteration & $\Delta^1$ & $\Delta^2$ & $\Delta^3$ & $\Delta^4$ & $\Delta^5$ \\
  \hline
  1&-0.203931& -0.205921& -0.203091& -0.205264& -0.201944\\
  \hline
  2&-0.105780& -0.102066& -0.103164& -0.102849& -0.108371\\
  \hline
  3&-0.109047& -0.105929& -0.105604& -0.105520& -0.111327\\
  \hline
  4&-0.108905& -0.105687& -0.105841& -0.105774& -0.111137\\
  \hline
  5&-0.108961& -0.105648& -0.105725& -0.105647& -0.111274\\
  \hline
\end{tabular}
 \caption{{\small Convergence of the delta for an American put basket option with $d=5$,
  $\rho=0.2$, $T=1$, $S_0 = 100$, $\sigma = 0.25$, $\delta=0.1$, $K=100$,
  $r=0.05$.}}
 \label{tab:delta_put_amer_5}
\end{table}

\subsection{Dupire's framework}
We consider a $d$-dimensional local volatility model in which the dynamics
under the risk-neutral measure of each asset is supposed to be given by 
\begin{equation*}
  dS_t^i = S_t^i ((r - \delta_i)dt + \sigma(t,S_t^i) dW^i_t) \qquad S_0 = (S_0^1, \dots, S_0^d)
\end{equation*}
where $W = (W^1, \dots, W^d)$ is defined and generated as in the Black-Scholes
framework. The local volatility function $\sigma$ we have chosen is of the
form
\begin{align}\label{eq:vol_locale}
  \sigma(t,x)=0.6(1.2-e^{-0.1t}e^{-0.001(xe^{rt}-s)^2})e^{-0.05 \sqrt{t}},
\end{align} with $s>0$. Since there exists a duality between the variables
$(t,x)$ and $(T,K)$ in Dupire's framework, one should choose $s$ equal to the
spot price of the underlying asset. Then, the bottom of the smile is located
at the forward money.  The parameters of the algorithm in this paragraph are
the following : $n=1000$, $M=30000$, $N=10$, $q=1$, $\eta=3$.

\paragraph{Pricing and Hedging European put basket options.}  We consider the
put basket option with payoff given by (\ref{eq:put_payoff}). The benchmark
price and delta are computed using the algorithm proposed
by~\cite{jourdainlelong}, which is based on Monte-Carlo methods.

\begin{figure}[h!t]
  \centering
  \includegraphics[scale=0.6]{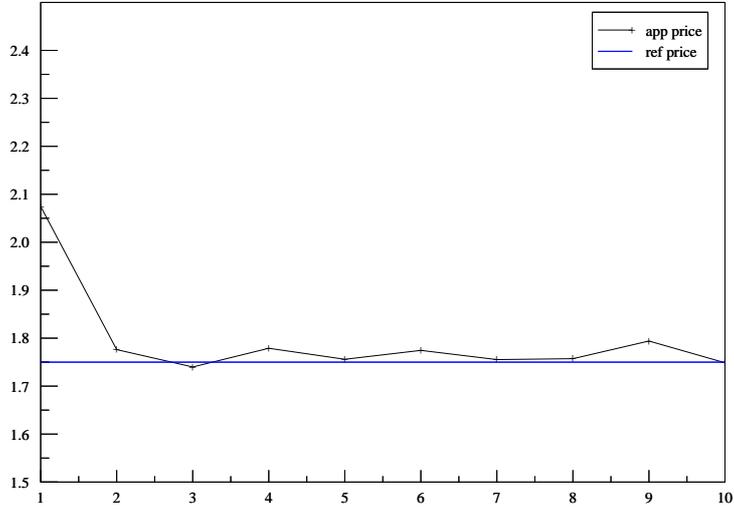}
  \caption{{\small Convergence of the algorithm for a European put basket option with
  $d=5$, $\rho=0$, $T=1$, $S_0 = 100$, $\delta=0$, $K=100$, $r=0.05$, $\eta=3$,
  $\omega=0$. The benchmark price computed with a high precision Monte--Carlo method
  yields $1.745899$ with a confidence interval of $(1.737899,1.753899)$.}}
  \label{fig:put_euro_vol_locale}
\end{figure}
Concerning the delta, we get at the last iteration the following vector
$\Delta=(-0.062403 -0.061271 -0.062437 -0.069120 -0.064743)$. The benchmark
delta is $-0.0625$.

\begin{figure}[h!t]
  \centering
  \includegraphics[scale=0.6]{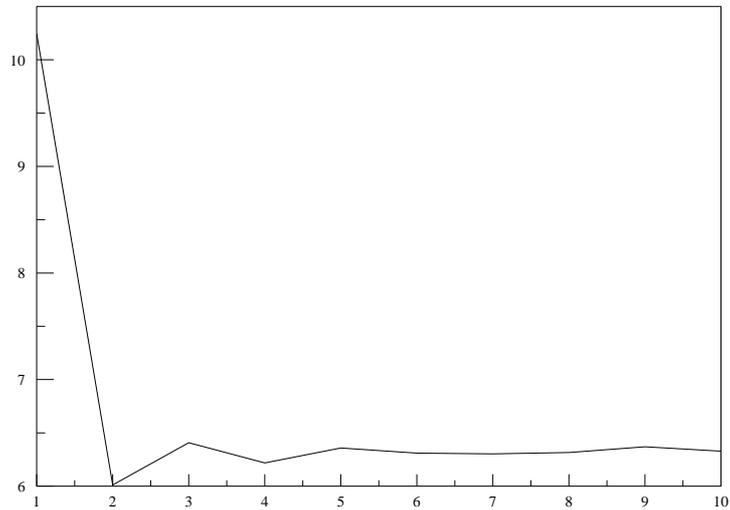}
  \caption{{\small Convergence of the algorithm for an American put basket option
  with $d=5$, $\rho=0$, $T=1$, $S_0 = 100$, $\delta=0.1$, $K=100$, $r=0.05$ and
  $\omega=1$.}}
  \label{fig:put_amer_vol_locale}
\end{figure}

\paragraph{Pricing and Hedging American put basket options.}
We still consider the put payoff given by (\ref{eq:put_payoff}). In the case
of American options in local volatility models, there is no
benchmark. However, Figure \ref {fig:put_amer_vol_locale} shows that the
algorithm converges after few iterations. We get a price around $6.30$. At the
last iteration, we get $\Delta=(-0.102159 -0.102893 -0.103237 -0.110546 -0.106442)$.

\subsection{Speed up.}

\begin{remark}
  Because in high dimensions, the sequential algorithm can run several hours before
  giving a price, we could not afford to run the sequential algorithm on the cluster
  to have a benchmark value for the reference CPU time used in the speed up
  measurements.  Instead we have computed speed ups as the ratio
  \begin{equation}
    \label{eq:speed-up-formula}
    \mbox{speed up} = \frac{\mbox{CPU time for $8$ processors / 8}}{\mbox{CPU time
    for $n$ processors $\times$ n}}
  \end{equation}
  This explains why we may get in the tables below some speed ups slightly greater
  than $1$.
\end{remark}

Our goal in this paper was to design a scalable algorithm for high dimensional
problems, so it is not surprising that the algorithm does not behave so well in
relatively small dimension as highlighted in Table~\ref{tab:speedup_put_amer_3}. In
dimension $3$, the speed ups are linear up to $28$ processors but then they
dramatically decrease toward zero: this can be explained by the small CPU load of the
computation of the correction term at a given point $(t_i, x_i)$. The cost of each
iteration of the loop line~\ref{loopi} of Algorithm~\ref{algodetail} is proportional
to $M p d^3$ and when $d$ is small so is $p$ --- the number of polynomials of total
degree less or equal than $\eta$. For instance, for $d=3$ and $\eta=3$, we have $p =
20$, which gives a small complexity for each iteration over $i$. Hence, when the
number of processors used increases, the amount of work to be done by its processor
between two synchronisation points decreases to such a point that most of the CPU
time is used for transmitting data or waiting. This explains why the speed ups
decrease so much. Actually, we were expecting such results as the parallel
implementation of the algorithm has been designed for high dimensional problems in
which the amount of work to be done by each processor cannot decrease so much unless
several dozens of thousands of processors are used. This phenomena can be observed in
Table~\ref{tab:speedup_put_amer_6} which shows very impressive speed ups: in
dimension $6$, even with $256$ processors the speed ups are still linear which
highlights the scalability of our implementation. Even though computational times may
look a little higher than with other algorithms, one should keep in mind that our algorithm not
only computes prices but also hedges, therefore the efficiency of the algorithm
remains quite impressive.

\begin{table}[ht!]
  \centering
  \begin{tabular}{ccc}
    \hline 
    Nb proc. &  Time    & Speed up \\ 
    \hline
    8        &  543.677 & 1        \\ 
    16       &  262.047 & 1.03737  \\ 
    18       &  233.082 & 1.03669  \\ 
    20       &  210.114 & 1.03501  \\ 
    24       &  177.235 & 1.02252  \\ 
    28       &  158.311 & 0.981206 \\ 
    32       &  140.858 & 0.964936 \\ 
    64       &  97.0629 & 0.70016  \\ 
    128      &  103.513 & 0.328267 \\ 
    256      &  162.936 & 0.104274 \\ 
    \hline
  \end{tabular}
  \caption{{\small Speed ups for the American put option with $d=3$, $r=0.02$, $T=1$,
  $\sigma=0.2$, $\rho=0$, $S_0 = 100$, $K=95$, $M=1000$, $N=10$, $K=10$, $n=2000$,
  $r=3$, $q=1$, $\omega=1$. See Equation~\eqref{eq:speed-up-formula} for the
  definition of the ``Speed up'' column.}}
  \label{tab:speedup_put_amer_3}
\end{table}

\begin{table}[ht!]
  \centering
  \begin{tabular}{ccc}
    \hline 
    Nb proc. &  Time    & Speed up \\ 
    \hline
    8        &  1196.79 & 1        \\ 
    16       &  562.888 & 1.06308  \\ 
    24       &  367.007 & 1.08698  \\ 
    32       &  272.403 & 1.09836  \\ 
    40       &  217.451 & 1.10075  \\ 
    48       &  181.263 & 1.10042  \\ 
    56       &  154.785 & 1.10457  \\ 
    64       &  135.979 & 1.10016  \\ 
    72       &  121.602 & 1.09354  \\ 
    80       &  109.217 & 1.09579  \\ 
    88       &  99.6925 & 1.09135  \\ 
    96       &  91.9594 & 1.08453  \\ 
    102      &  85.6052 & 1.0965   \\ 
    110      &  80.2032 & 1.08523  \\ 
    116      &  75.9477 & 1.08676  \\ 
    128      &  68.6815 & 1.08908  \\ 
    256      &  35.9239 & 1.04108  \\ 
    \hline
  \end{tabular}
  \caption{{\small Speed ups for the American put option with $d=6$, $r=0.02$, $T=1$,
  $\sigma=0.2$, $\rho=0$, $S_0 = 100$, $K=95$, $M=5000$, $N=10$, $K=10$, $n=2000$,
  $r=3$, $q=1$, $\omega=1$. See Equation~\eqref{eq:speed-up-formula} for the
  definition of the ``Speed up`` column.}}
  \label{tab:speedup_put_amer_6}
\end{table}

\section{Conclusion}

In this work, we have presented a parallel algorithm for solving BSDE and applied it
to the pricing and hedging of American option which remains a computationally
demanding problem for which very few scalable implementations have been studied. Our
parallel algorithm shows an impressive scalability in high dimensions. To improve the
efficiency of the algorithm, we could try to refactor the extrapolation step to make
it more accurate and less sensitive to the curse of dimensionality.

\clearpage
\bibliographystyle{abbrvnat}
\bibliography{biblio}

\end{document}